\documentclass[12pt,a4paper]{article}

\usepackage{amsfonts}
\usepackage{amsmath}
\usepackage{amsthm}
\usepackage{amssymb}

\theoremstyle{definition} \newtheorem{defn}{Definition}[section]
\theoremstyle{plain} \newtheorem{thm}[defn]{Theorem}
\theoremstyle{plain} \newtheorem{propn}[defn]{Proposition}
\theoremstyle{plain} \newtheorem{lemma}[defn]{Lemma}
\theoremstyle{plain} \newtheorem{cor}[defn]{Corollary}
\theoremstyle{plain} \newtheorem*{claim}{Claim}

\theoremstyle{definition}  \newtheorem*{ex}{Example}
\theoremstyle{definition}  
\theoremstyle{remark}  \newtheorem*{rmk}{Remark}

\newcommand {\R} {\mathbb{R}}
\newcommand {\Q} {\mathbb{Q}}

\newcommand {\N} {\mathbb{N}}
\newcommand {\Rbar} {\bar{\mathbb{R}}}
\newcommand {\cl} [1] {\overline{#1}} 

\newcommand {\fr} {\textrm{fr}}
\newcommand {\la} {\langle}
\newcommand {\ra} {\rangle}

\newcommand {\abar} {\bar{a}}
\newcommand {\bbar} {\bar{b}}

\newcommand {\xbar} {\bar{x}}
\newcommand {\ybar} {\bar{y}}

\newcommand {\etabar} {\bar{\eta}}
\newcommand {\ytild} {\tilde{y}}
\newcommand {\atild} {\tilde{a}}

\newcommand {\C} {\mathcal{C}}

\renewcommand {\epsilon}{\varepsilon}

\newcommand {\M} {\mathcal{M}}

\newcommand {\curlyN} {\mathcal{N}}
\newcommand {\Mbar} {\bar{M}}

\newcommand {\F} {\mathcal{F}}
\newcommand {\G} {\mathcal{G}}

\renewcommand {\t} [1] {\widetilde{#1}} 

\newcommand {\grad} {\nabla}
\renewcommand {\d}[2] {\frac{\partial #1}{\partial #2}}

\title{Locally Polynomially Bounded Structures}
\author{G. O. Jones \and A. J. Wilkie}
\begin{document}
\maketitle
In his proof of the model completeness of the real exponential field
(\cite{rexp}), the second author develops a theory of Noetherian differential
rings of definable functions, and studies varieties defined by these
functions. One of the main results of this theory is Theorem 5.1 which provides a
method for constructing points on such varieties. 

Our aim in this paper is to prove a version of this theorem without
the Noetherianity assumption. Instead we suppose that the functions
considered are definable in an expansion of a real closed field, $\M$
say, which is definably complete (see \cite{millerIVP},\cite{tamara}) and further, that
the functions are what we call \emph{locally tame}. We will give
precise definitions later, but the idea is that certain
restrictions (to bounded boxes) of the (total) functions considered, are
definable in a fixed o-minimal polynomially bounded reduct of $\M$. Then we
can use Miller's results (\cite{millerqa},\cite{millerpowerbdd}) to bound orders of
vanishing and it is this that makes up for the lack of Noetherianity.

After proving the main result, we specialize to the o-minimal
situation. We call an o-minimal structure $\M$ with model complete
theory \emph{locally polynomially bounded} if the reduct generated by
all restrictions of the basic functions to bounded open boxes is
polynomially bounded. We show that being locally polynomially bounded
is preserved under elementary equivalence. Combining this with model completeness and the main
theorem, we show that definable functions are piecewise implicitly
defined over the basic functions in the language. This implies that
these structures have smooth cell decomposition. Under a further
assumption on these basic functions, this gives uniform control over
the derivatives of definable functions.
\section{Locally tame functions}\label{lpb}
Let $\bar{M} = \langle M, <, +, \cdot, 0, 1 \rangle$ be a fixed real closed
field and let $\mathcal{M} = \langle \bar{M}, \ldots \rangle$ be an arbitrary,
 but fixed, expansion of $\mathcal{M}$. We also fix an o-minimal,
polynomially bounded reduct $\mathcal{M}_0$ of $\mathcal{M}$ such that
$\mathcal{M}_0$ is also an expansion of $\bar{M}$. 
We use \emph{definable} to
mean definable with parameters and \emph{0-definable}
to mean definable  without parameters, and unless we specifically mention
 another structure, we are referring to definability in $\mathcal{M}$.
\begin{defn}Suppose that $f:U \to M$ is a definable function
on some open $U \subseteq M^n$. We say that $f$ is \emph{locally tame}
if $f$ is smooth (i.e. infinitely differentiable on $U$ in the sense of the
usual $\epsilon-\delta$ definition formulated in $\mathcal{M}$) and, for
every open box $B \subseteq M^n$ having sides
of length $\leq 1$ and satisfying $\bar{B} \subseteq U$,
 we have that $f|_B$ is definable in $\mathcal{M}_0$.
\end{defn}
\begin{ex}
Suppose that $\mathcal{M} = \langle \bar{\mathbb{R}}, \exp
\rangle$ and $\mathcal{M}_0 = \langle \bar{\mathbb{R}}, \exp| [0, 1] \rangle$.
Then $\exp$ is locally tame. Now consider the function
\begin{eqnarray*}
g:\R&\to&\R \\
t&\mapsto&
\begin{cases}
\exp (-1/t^2) & t\ne 0\\
0 & t=0.
\end{cases}
\end{eqnarray*}
This function is smooth and definable, but, by the following result,
it is not locally tame.
\end{ex}
\begin{propn}\label{noflatpoints}
 Let $f:U \to M$ be a locally tame function. Then the set of flat
 points of $f$ (i.e. the points at which all derivatives of $f$ of all
 orders vanish) is definable and is both open and closed in
 $U$. Further, if $B \subseteq U$ is any open box having sides of length at
 most $1$, and if $B$ contains a flat point of $f$, then $f$ vanishes
 throughout $B$.
\end{propn}
\begin{proof}
Let $X$ be the set of all flat points of $f$ and let 
$$
Y=\{ \xbar \in U: \textrm{ there is an open box around } \xbar \textrm{ on
  which }f\textrm{ vanishes}\}.
$$
Clearly we have $Y\subseteq X$. Suppose that $\abar \in X$ and let $B$
be a box containing $\abar$ with sides of length $\leq 1$ such that
$\cl{B}\subseteq U$. Since $f$ is locally tame, the restriction $f|_B$
is definable in $\M_0$. This structure is o-minimal and polynomially
bounded so a result of Miller's (\cite{millerqa}) shows that $\abar
\in Y$. So $Y=X$ and $X$ is definable and since $Y$ is open, so is
$X$.

Now suppose for a contradiction that $X$ is not closed in $U$. Then
there is some $\bbar\in U$ such that $\bbar\in \fr X$. Fix $\alpha \in
\N^n$. There are points arbitrarily close to $\bbar$ at which $f$ is
flat. At such a point, $\xbar$ say, $D^{\alpha}f(\xbar)=0$. Hence
$D^{\alpha}f(\bbar)=0$. Since $\alpha\in\N^n$ was arbitrary, it
follows that $f$ is flat at $\bbar$. So $\bbar \in X$ which is a
contradiction.  

Finally, if
$B \subseteq U$ is an open box having sides of  length at most $1$ and
containing a flat point of $f$, then we may apply the above argument in
the o-minimal, polynomially bounded structure $\mathcal{M}_0$ to the
function $f|B$ and use the fact that $B$ is definably connected.
\end{proof}
We now suppose that we have, for each $n\ge 1$, a $\Q$-algebra $R_n$
of locally tame functions $f:M^n\to M$, which is closed under partial
differentiation. We will also assume that $R_n\subseteq R_{n+1}$ (in
the obvious sense) and that
$$
\Q[X_1,\ldots,X_n]\subseteq R_n.
$$
Before giving our main result, we recall some notation from
\cite{rexp}. Let $f\in R_n$. We define $\grad f:M^n\to M^n$ by
\begin{equation*}
\grad f(\abar) :=\la \d{f}{x_1}(\abar),\ldots,\d{f}{x_n}(\abar)\ra\qquad
\textrm{for } \abar \in M^n.
\end{equation*}
Note that $\grad f \in R_n^n$. For $p\ge 1$ and $f_1,\ldots,f_p \in R_n$ we let
$$
V_n(f_1,\ldots,f_p):=\{\xbar \in M^n:f_1(\xbar)=\cdots=f_p(\xbar)=0\}
$$
and 
$$
V^{reg}_n(f_1,\ldots,f_p):=\{ \xbar \in V_n(f_1,\ldots,f_p):\grad
f_1(\xbar),\ldots,\grad f_p(\xbar) \textrm{ are linearly independent}\}.
$$
Here, linear independence is in the $\Mbar$ vector space $M^n$. The
\emph{Jacobian matrix} of $f_1,\ldots,f_p$ is the matrix
$$
J_n(f_1,\ldots,f_p):= \left( \begin{array}{l} 
\grad f_1 \\
\vdots \\
\grad f_p \end{array} \right).
$$
The rows of $J_n(f_1,\ldots,f_p)$ are linearly independent when
evaluated at $\abar \in M^n$ if and only if $p\le n$ and there is a
$p\times p$ submatrix whose determinant is non-zero when evaluated at
$\abar$. So, if we let $Q=Q_{n,f_1,\ldots,f_p}\in R_n$ be the sum of
squares of all such determinants we have
\begin{equation}\label{detwitness}
\textrm{for all }\abar\in M^n, \abar \in
V^{reg}_n(f_1,\ldots,f_p)\leftrightarrow \abar \in V_n(f_1,\ldots,f_p)
\textrm{ and } Q(\abar)>0.
\end{equation}
\begin{lemma}\label{Vclosed}
 If we regard $f_1,\ldots,f_p$ as elements of $R_{n+1}$, then there is
 a function $f_{p+1}\in R_{n+1}$ such that 
$$
V_{n+1}(f_1,\ldots,f_{p+1})=V_{n+1}^{reg}(f_1,\ldots,f_{p+1}),
$$
 and $V_{n+1}(f_1,\ldots,f_{p+1})$ projects onto
 $V^{reg}_n(f_1,\ldots,f_p)$. In particular, $V_{n+1}^{reg}(f_1,\ldots,f_{p+1})$ is closed in $M^{n+1}$.
\end{lemma}
\begin{proof} Let $f_{p+1}(x_1,\ldots,x_{n+1})=x_{n+1}\cdot
  Q(x_1,\ldots,x_n)-1$ where $Q$ is the function defined before the
  Lemma. Then by (\ref{detwitness}),
$$
\la \abar,a_{n+1}\ra \in V^{reg}_{n+1}(f_1,\ldots,f_{p+1})\leftrightarrow
\abar \in V_n^{reg}(f_1,\ldots,f_p) \textrm{ and } a_{n+1}=Q(\abar)^{-1}.
$$
An easy calculation shows that for such $\la \abar,a_{n+1}\ra$, we have
$Q_0(\abar,a_{n+1})\ge Q(\abar)^3$, where
$Q_0:=Q_{n+1,f_1,\ldots,f_{p+1}}$. So  
$$
V_{n+1}(f_1,\ldots,f_{p+1})=V_{n+1}^{reg}(f_1,\ldots,f_{p+1}).
$$
as required.
\end{proof}
\begin{thm} \label{5point1}Assume that $\mathcal{M}$ is
definably complete (i.e. every definable subset of $M$ with an upper bound
has a least upper bound).
Suppose that $n\ge 1 $ and that $f\in R_n$ is such that $V(f)$ is
  nonempty. Then there exist $m\ge 0$ and $f_1,\ldots,f_{n+m}\in  R_{n+m}$
  such that
$$
V_{n+m}^{reg}(f_1,\ldots,f_{n+m})\cap V_{n+m}(f)\ne \emptyset.
$$
Here we regard $f$ as an element of $R_{n+m}$, so $V_{n+m}(f)=V_n(f)\times M^m$.
\end{thm}
\begin{proof} If $f$ vanishes identically then we let
  $m=0$ and $f_i(x_1,\ldots,x_n)=x_i$ so that we have 
  $V^{reg}_n(f_1,\ldots,f_n)\cap V_n(f)=\{ \cl{0}\}$. So we may
  suppose that $f$ is not identically zero.

We will show by induction on $p$, for $1\le p\le n$, that 
\begin{equation}\label{starp}
\begin{array}{l}
\textrm{ there exist }m\ge 0 \textrm{ and } f_1,\ldots,f_{p+m}\in R_{p+m} \textrm{
  such that } \\
 V_{n+m}^{reg}(f_1,\ldots,f_{p+m})\cap V_{n+m}(f)\ne
\emptyset.
\end{array}
\end{equation}
Suppose first that $p=1$. We choose any point $\abar \in
V_n(f)$. Since $\M$ is definably complete, the set $M$ is definably
connected (see \cite{millerIVP}) and a simple argument shows that
$M^n$ is also definably connected. Hence since $f$ is locally tame and is not
identically zero, Proposition \ref{noflatpoints} gives an $\alpha \in
\N^n$ such that, with $f_1=D^{\alpha}f$, we have $\abar \in
V^{reg}_n(f_1)$. This proves \eqref{starp} for $p=1$, with $m=0$.

Now suppose that $p$ is such that $1\le p <n$ and that \eqref{starp}
holds for $p$. Then we have $m\ge 0$ and $f_1,\ldots,f_{p+m}\in R_{n+m}$
such that
$$
V_{n+m}^{reg}(f_1,\ldots,f_{p+m})\cap V_{n+m}(f)\ne \emptyset.
$$
\emph{Case 1.} There is some $\abar \in
  V_{n+m}^{reg}(f_1,\ldots,f_{p+m})\cap V_{n+m}(f)$ such that $f$ is
  not identically zero on $B\cap V^{reg}_{n+m}(f_1,\ldots,f_{p+m})$
  for any open box $B\subseteq M^{n+m}$ with $\abar\in B$.

Since $\abar \in V^{reg}_{n+m}(f_1,\ldots,f_{p+m})$, there is some
$(p+m)\times (p+m)$ submatrix of $J_{n+m}(f_1,\ldots,f_{p+m})$ whose
determinant is non-zero at $\abar$. We will assume that this submatrix
consists of the last $(p+m)$ columns, and write $\Delta$ for its
determinant. Note that $\Delta$ is a function in $R_{n+m}$. For
$\ybar=\la y_1,\ldots,y_{n+m}\ra \in M^{n+m}$, we let $\ytild:=\la
y_1,\ldots,y_{n-p}\ra$. Since the functions $f_1,\ldots,f_{p+m}$ are
locally tame, there is an open box $B_0\subseteq M^{n+m}$ such that
$\abar \in B_0$ and $f_1|_{B_0},\ldots,f_{p+m}|_{B_0}$ are definable
in $\M_0$. By the implicit function theorem, applied
in the o-minimal structure $\M_0$ (see \cite{thebook},
Chapter 7) there is an open box $U\subseteq M^{n-p}$ with $\atild\in
U$ and a smooth map $\phi:U\to M^{p+m}$, definable in
$\M_0$, such that 
\begin{enumerate}
\item[(i)] $\phi(\atild)=\la a_{n-p+1},\ldots,a_{n+m}\ra$,
\item[(ii)] $\{ \la \ytild,\phi(\ytild)\ra:\ytild \in U\} =B\cap V^{reg}_{n+m}(f_1,\ldots,f_{p+m})$
\end{enumerate}
for some open box $B\subseteq M^{n+m}$ with $\abar \in B$. We may
suppose that $\Delta$ has no zeroes in $B$. Since $f$ is locally tame,
$f|_B$ is definable in $\M_0$ and hence so is the
function
\begin{eqnarray*}
g:U&\to &M \\
\ytild& \mapsto &f(\ytild,\phi(\ytild )).
\end{eqnarray*}
 Now, by the hypothesis of case 1 and (i) and (ii) above, $g$ is not
identically zero on $U$, and as $\M_0$ is polynomially bounded, there
is some $\alpha \in \N^{n-p}$ such that $g^*:=D^{\alpha}g$ vanishes at
$\atild$ but, for some $j=1,\ldots,{n-p},\d{g^*}{y_j}$ does not.

Now we have 
\begin{eqnarray*}
f_i(\ytild,\phi(\ytild ))&=&0  \textrm{ for }i=1,\ldots,m+p \textrm{
  and }\ytild \in U,\\
g(\ytild)&=&f(\ytild,\phi(\ytild ))  \textrm{ for } \ytild \in U,
\end{eqnarray*}
and by differentiating these relations, we obtain a function $F\in
R_{n+m}$ such that
\begin{equation*}
g^*(\ytild )=\frac{F(\ytild, \phi (\ytild ))}{\Delta (\ytild,\phi(
  \ytild))^d}  \textrm{ for all } \ytild \in U
\end{equation*}
for some $d$. We also have that $F(\atild,\phi(\atild ))=F(\abar)=0$,
and since $\d{g^*}{y_j}(\atild )\ne 0$, it follows from Lemma 4.7 in
\cite{rexp} that $\grad f_1(\abar),\ldots,\grad f_p (\abar),\grad
F(\abar)$ are linearly independent. So we obtain (\ref{starp}) for
$p+1$ by taking $f_{p+m+1}=F$ and not changing $m$.

\emph{Case 2.} Not case 1.

By Lemma \ref{Vclosed} we may suppose (after increasing $m$) that 
$V^{reg}_{n+m}(f_1, \ldots, f_{p+m}) = V_{n+m}(f_1, \ldots, f_{p+m})$. Let
$\mathcal{C} = 
V^{reg}_{n+m}(f_1, \ldots, f_{p+m}) \cap V_{n+m}(f)$. Then $\mathcal{C}$ is
nonempty (by (2)) and closed in $M^{n+m}$. Now, if we can find some
$h \in R_{n+m}$ which has a zero in $\mathcal{C}$ but is not identically
zero on $B \cap V^{reg}_{n+m}(f_1, \ldots, f_{p+m})$ for any open box $B$
containing this zero, then we can apply the method of
Case 1 to $h$ and we will be done.

To find such an $h$ we proceed as in the proof of Theorem 5.1 in
\cite{rexp}. Let $\etabar=\la \eta_1,\ldots,\eta_{n+m}\ra \in
\Q^{n+m}$. Then, since $\C$ is closed, there is a point $\bbar\in
\C$ at minimum distance from $\etabar$. (This follows
easily from the definable completeness of $\mathcal{M}$.) Let
$H_{\bar{\eta}}(\bar{x}) := \Sigma(x_i - \eta_i)^2$. 
Then $H_{\bar{\eta}} \in R_{n+m}$ and the function
$H_{\bar{\eta}}|\mathcal{C}$ has a minimum at $\bar{b}$. 
However, by the hypotheses
of Case 2, $\mathcal{C}$ coincides with  $V^{reg}_{n+m}(f_1, \ldots, f_{p+m})$
on some open box in $M^{n+m}$ containing the point $\bar{b}$ and hence, by
the method of Lagrange multipliers (see 4.10 in \cite{rexp};
we should also remark that we may work in the o-minimal structure
$\mathcal{M}_0$ at this point), the vectors
$\nabla f_1(\bar{b}), \ldots,
\nabla f_{p+m}(\bar{b}), \nabla H_{\bar{\eta}} (\bar{b})$ 
are linearly dependent. Now, by (1), this is equivalent to the vanishing
at $\bar{b}$ of the function $Q_{\bar{\eta}} := Q_{n+m, f_1 , \ldots ,
f_{p+m} , H_{\bar{\eta}}} \in R_{n+m}$. Now consider the function 
$\tilde{f} := Q_{\bar{\eta}}^2
+ f^2$. Either it will serve as the required function $h$, or else it too
satisfies the same hypothesis of Case 2 as did $f$  
(including the fact that $V^{reg}_{n+m}(f_1, \ldots, f_{p+m}) \cap V_{n+m}
(\tilde{f}) \neq \emptyset$). 

So by successive repetition of this argument
we either succeed in finding a suitable $h$, or else for any positive
integer $r$ and any sequence of
points $\bar{\eta}_{1}, \ldots \bar{\eta}_{r} \in 
\mathbb{Q}^{n+m}$, we find a point $\bar{c} \in \mathcal{C}$ such that
for each $i = 1, \ldots ,r$, the vector
$\nabla H_{\bar{\eta}_{i}}(\bar{c})$ lies in the vector space spanned by
$\nabla f_1(\bar{c}), \ldots, \nabla f_{p+m}(\bar{c})$. However, for
$\bar{\eta} \in \mathbb{Q}^{n+m}$ one calculates that
$\nabla H_{\bar{\eta}}(\bar{c}) 
= \langle 2(c_1 - \eta_{1}), \ldots ,2(c_{n+m} -  
\eta_{n+m}) \rangle$. Thus, if we take $r = n+m+1$, $\bar{\eta}_1 = \bar{0}$
and $\bar{\eta}_2 , \ldots ,  \bar{\eta}_{n+m+1}$ to be any basis for
$\mathbb{Q}^{n+m}$ we see that (for any $\bar{c} \in M^{n+m}$), the set
$\{\nabla H_{\bar{\eta}_1}(\bar{c}), \ldots , 
\nabla H_{\bar{\eta}_{n+m+1}}(\bar{c})\}$ spans $M^{n+m}$, contradicting
the fact that $p+m < n+m$. Thus we will find a suitable $h$ 
and this completes the proof of Theorem \ref{5point1}.
\end{proof}
Examples of definably complete structures include any structure
elementarily equivalent to an expansion of $\bar{\mathbb{R}}$, and any  
o-minimal structure. See \cite{millerIVP} and \cite{tamara} for discussions of definably complete 
structures.
\section{Locally polynomially bounded structures}\label{lpbstructures}
From now on in this paper we consider the case that $\mathcal{M} = 
\langle \bar{M}, \mathcal{F} \rangle$ and $\mathcal{M}_0 =
\langle \bar{M}, \mathcal{F}^{res} \rangle$, where $\mathcal{F}$ is a
family of smooth functions $f:M^n \to M$, for various $n$, and where
$\mathcal{F}^{res}$ denotes the collection of all functions of the form
$f|_B$ for $f \in \mathcal{F}$ and $B$ an open box in $M^n$. We call $\mathcal{M}$ \emph{locally polynomially bounded}
 (LPB) if $\mathcal{M}$ is o-minimal and has a model complete theory (as
well as $\mathcal{M}_0$ being polynomially bounded). 

Thus, the example
of the real exponential field discussed in the previous section is LPB. More
generally suppose that $\tilde{\mathbb{R}}$ is a polynomially bounded
o-minimal expansion of $\Rbar$ with smooth cell decomposition and
  that the restricted exponential function, $\exp|_{[0,1]}$, is definable in $\t{\R}$. Let $\F$
  denote the collection of all total smooth definable
  functions. By a Theorem of van den Dries and Speissegger (Theorem B
  in \cite{vdDSp}) the structure $\la \Rbar,\F,\exp\ra$ is model
  complete and hence LPB.

Now, let $\mathcal{M}$ be an LPB
structure with $\mathcal{M}_0$ and $\mathcal{F}$  as described above. Let
$\mathcal{N} = \langle \bar{N}, \mathcal{G} \rangle$ be a structure for the 
same language (with $\bar{N}$ a real closed field) 
and form $\mathcal{G}^{res}$ and $\mathcal{N}_0$ in the analogous way.

\begin{thm} \label{elemequiv} If $\mathcal{M} \equiv \mathcal{N}$ then $\mathcal{N}$
is also LPB.
\end{thm}
\begin{proof} 
We must show that $\curlyN_0$ is polynomially bounded. We will first
show it is power
bounded, in the sense of \cite{millerpowerbdd}, so suppose that it is
not. Then by
Miller's Dichotomy Theorem (\cite{millerpowerbdd}), there is an
exponential function $E:N\to N$ which is 
$0$-definable in $\curlyN_0$. This means that there is some formula, $\Phi
(F_1,\ldots,F_n,x,y)$ say, in the language of ordered rings together
with $n$ function variables (of various arities) but only first order quantifiers,
functions $g_1,\ldots,g_n\in \G$ (of the corresponding arities) and bounded open boxes
$B_1,\ldots,B_n$ (in the corresponding spaces) such that 
$$
 \textrm{ for all }a,b\in N, E(a)=b \textrm{ if and only if } \curlyN \models \Phi
(g_1|_{B_1},\ldots,g_n|_{B_n},a,b).
$$
We now write $\Phi(x,y)$ for $\Phi(F_1,\ldots,F_n,x,y)$ and let  $\Psi
(F_1,\ldots,F_n)$ be the formula
\begin{eqnarray*}
\forall x\exists !y \Phi(x,y)\land \Phi(0,1)\land \forall x,x',y,y'
\big[ (x<x' \land \Phi(x,y)\land\Phi(x',y'))\to \\
(y<y'\land \Phi(x+x',y\cdot y'))\big].
\end{eqnarray*}
Then 
$$
\curlyN \models \exists B_1,\ldots,B_n \Psi (g_1|_{B_1},\ldots,g_n|_{B_n}).
$$
Now quantification over boxes is first order, as we can quantify over
the corners. 
 So, by the elementary equivalence of $\M$ and $\curlyN$, we have
$$
\M \models \exists B_1,\ldots,B_n \Psi (f_1|_{B_1},\ldots,f_n|_{B_n})
$$
where the $f_1,...,f_n 
\in \mathcal{F}$ correspond to $g_1,...,g_n \in \mathcal{G}$.
Hence an exponential function is definable in $\M_0$, contradicting the fact that $\M$ is locally
polynomially bounded.

So $\curlyN_0$ is power bounded. We now need to show that it is
polynomially bounded. 
\begin{claim}Suppose that for any formula
  $\Phi(F_1,\ldots,F_n,x,y)$ (in the language of ordered rings, together
with $n$ function variables but only first order quantifiers)
  and any collection of functions $f_1,\ldots,f_n \in \F$, the formula
  $\Phi(f_1|_{B_1},\ldots,f_n|_{B_n},x,y)$ defines in $\M$ the graphs of only
  finitely many power functions as the boxes $B_1,\ldots,B_n$
  vary. Then $\curlyN_0$ is polynomially bounded.
\end{claim}
\begin{proof}
Note that it suffices to show that there is no non-polynomially
bounded power function definable without parameters in $\curlyN_0$.
So, suppose that $g_1\ldots,g_n\in\G$ and $B_1,\ldots,B_n$ are open boxes
such that the formula $\Phi(g_1|_{B_1},\ldots,g_n|_{B_n},x,y)$ defines a
power function, $x^{\alpha}$ say, in $\curlyN_0$. By the hypothesis of the Claim and the fact that
$\M_0$ is polynomially bounded, there is a $k\in \N$ such that
the sentence
\begin{eqnarray*}
\forall B_1,\ldots,B_n (\textrm{if
}\Phi(f_1|_{B_1},\ldots,f_n|_{B_n},x,y)\textrm{ defines a power function }\\
 \textrm{then this function is bounded by }x^k)
\end{eqnarray*}
holds in $\M$, where the $ f_1,...,f_n 
\in \mathcal{F}$ correspond to $ g_1,...,g_n \in \mathcal{G}$. (To see that the set of boxes for which
$\Phi(f_1|_{B_1},\ldots,f_n|_{B_n},x,y)$ defines a power function is
definable, write out a formula analogous to the formula $\Psi$
above.) Hence this sentence is true in $\curlyN$ and so $\alpha\le k$
and $\curlyN_0$ is polynomially bounded.
\end{proof}

We will now establish the hypothesis of the Claim, so fix a formula
$\Phi(F_1,\ldots,F_n,x,y)$ and functions $f_1,\ldots,f_n\in\F$. Let $K$ be the (definable)
set of exponents of power functions defined by
$\Phi(f_1|_{B_1},\ldots,f_n|_{B_n},x,y)$ as the boxes $B_1,\ldots,B_n$
vary and suppose for a contradiction that $K$ contains a nonempty open
interval, $J$ say. Using definable choice and monotonicity (in the
o-minimal structure $\M$) there is a
bounded subinterval $J_0$ say, with $\cl{J}_0\subseteq J$ and a
continuous definable function $G$ on $\cl{J}_0$, whose values are
$n$-tuples of boxes, such that 
\begin{equation}\label{powerfns}
\begin{array}{l}
\textrm{for all }\alpha \in \cl{J}_0, G(\alpha)=\la B_1^{\alpha},\ldots,B_n^{\alpha}\ra \textrm{ is such that }\\
\Phi(f_1|_{B_1^{\alpha}},\ldots,f_n|_{B_n^{\alpha}},x,y) \textrm{
  defines } y=x^{\alpha} \textrm{ for }x>0.
\end{array}
\end{equation}
Since $\cl{J}_0$ is a closed bounded interval, $G$ is bounded and we
may take bounded open boxes $D_1,\ldots,D_n$ such that
$B_i^{\alpha}\subseteq D_i$, for all $\alpha \in \cl{J_0}$ and
$i=1,\ldots,n$. Now, if we repeat the above argument with the
structure $\la \Mbar,f_1|_{D_1},\ldots,f_n|_{D_n}\ra$ in place of
$\M$, we obtain an interval $J_1$ and a function $G_1$ defined in the structure $\la
\Mbar,f_1|_{D_1},\ldots,f_n|_{D_n}\ra$ such that (\ref{powerfns})
holds with $J_1,G_1$ in place of $J_0,G$. Hence the function $\la
x,y\ra \mapsto x^y$ with $x>0$ and $y\in J_1$ is definable in $\la
\Mbar,f_1|_{D_1},\ldots,f_n|_{D_n}\ra$ and hence in $\M_0$. But this
is impossible, by the proof of 4.2 in \cite{millerpowerbdd}, as
$\M_0$ is polynomially bounded.
\end{proof}
\section{Consequences of model completeness}
For the remainder of the paper, we fix an LPB structure, $\M=\la
\Mbar,\F\ra$. Let $\t{\F}$ be the smallest collection of functions containing $\F$
and all polynomials over $\Q$ and closed under the $\Q$-algebra
operations and under partial differentiation. For each $n\ge 1$, let
$R_n$ be the $\Q$-algebra consisting of all $n$-ary functions in
$\t{\F}$. Then each $R_n$ is closed under partial differentiation and
consists of locally tame functions, so the results of the first
section apply. These results also apply to the rings $R^{\abar}_n$,
for $\abar\in M^p$, consisting of all functions of the form
$\xbar\mapsto f(\abar,\xbar)$ for some $f\in R_{p+n}$.

\begin{defn}Let $\abar\in M^p$ and $b\in M$. We say that $b$ is
$\F$-\emph{defined over} $\abar$ if there exist $n\ge
1,f_1,\ldots,f_n\in R^{\abar}_n$ and $b_1,\ldots,b_n\in M$ with
$b=b_i$ for some $i$, such that
$$
\bbar\in V^{reg}_n (f_1,\ldots,f_n).
$$
\end{defn}
The following is an easy consequence of the model completeness of $\la
\Mbar,\F\ra$ and \ref{5point1}, together with a standard trick on
representing definable sets as projections of zero sets.
\begin{thm} Let $\abar\in M^p,b\in M$. Then $b$ is in the definable
  closure of $\abar$ if and only if $b$ is $\F$-defined over
  $\abar$. In particular, ``$\F$-defined over'' is a pregeometry.
\end{thm}
\begin{defn}We say that a $0$-definable function $f:U\to
  M$, where $U\subseteq M^n$ is open, is \emph{implicitly defined
    over} $\F$  if there exist $m\ge 1$, functions $g_1,\ldots,g_m\in
  R_{n+m}$ and $0$-definable functions $\phi_1,\ldots,\phi_m:U\to M$
  such that
\begin{enumerate}
\item[(1)] $f=\phi_i$, for some $i=1,\ldots,m$,
\item[(2)] $\la \phi_1(\xbar),\ldots,\phi_m(\xbar)\ra \in
  V^{reg}_n(g_1(\xbar,\cdot),\ldots,g_m(\xbar,\cdot))$, for all $\xbar
  \in U$.
\end{enumerate}
\end{defn}
\begin{cor}\label{piecewiseimplicit}
 If $\abar \in M^n$ is generic (for the pregeometry given
  by definable closure) and $f:U\to M$ is a $0$-definable function on
  a neighbourhood $U$ of $\abar$, then there is an open $0$-definable
  $V\subseteq U$ with $\abar\in V$, such that $f|_V$ is implicitly defined over $\F$.
\end{cor}
\begin{proof} Since $f(\abar)$ is in the definable closure of $\abar$,
  it follows from the previous Theorem that there exist $m\ge 1$,
  functions $g_1,\ldots,g_m\in R_{n+m}$ and a tuple $\la
  b_1,\ldots,b_m\ra \in M^m$ such that $f(\abar)\in \{
  b_1,\ldots,b_m\}$ and
$$
\bbar\in V^{reg}_m(g_1(\abar,\cdot),\ldots,g_m(\abar,\cdot)).
$$
Consider the $0$-definable set
$$
X:= \{ \la \xbar,\ybar\ra \in M^{n+m}: \ybar\in
V^{reg}_m(g_1(\xbar,\cdot),\ldots,g_m(\xbar,\cdot))\}.
$$
 For each $\xbar$ there are at most finitely many $\ybar$
such that $\la \xbar,\ybar\ra\in X$. Hence by cell decomposition and
the fact that $\abar$ is generic there is an open cell, $C$ say,
containing $\abar$ and $0$-definable functions $\phi_1,\ldots,\phi_m:C\to M$ such that 
$$
\la \phi_1(\xbar),\ldots,\phi_m(\xbar)\ra \in
V^{reg}_m(g_1(\xbar,\cdot),\ldots,g_m(\xbar,\cdot))
$$
and $\phi_i(\abar)=f(\abar)$ for some $i$. Then as $\abar$ is
generic, $\phi_i$ and $f$ agree on some open neighbourhood $V$ of
$\abar$ and so $f|_V$ is implicitly defined over $\F$.
\end{proof}
Using this Corollary and Theorem \ref{elemequiv}, a standard
compactness argument yields the following:
\begin{cor} Suppose that $f:U\to M$ is a $0$-definable function on an
  open set $U\subseteq M^n$. Then there are $0$-definable open sets
  $U_1,\ldots,U_k\subseteq U$ with $\dim (U\setminus \bigcup_{i=1}^k
  U_i)<n$ such that $f|_{U_i}$ is implicitly defined over $\F$, for
  each $i$.
\end{cor}
Now the implicit function theorem implies that functions which are
implicitly defined over $\F$ are smooth, and so we have:
\begin{cor}
Locally polynomially bounded structures have smooth cell
decomposition.
\end{cor}
\section{Controlling the derivatives}
\begin{defn}
A smooth definable function $f:U\to M$ on an
open set $U\subseteq M^n$ is said to have \emph{controlled
  derivatives} if there exists a definable continuous function $\omega:U\to
M_{\ge 0}$ and $C_i\in M,E_i\in \N$, for each $i\in\N$ such that
$$
|D^{\alpha}f(\xbar)|\le C_{|\alpha|}\cdot \omega
(\xbar)^{E_{|\alpha|}} \textrm{ for all }\alpha\in \N^n\textrm{ and
}\xbar \in U.
$$
We say that such an $\omega$ is a \emph{control function} for $f$ and
that $\{ \omega,C_i,E_i\}$ is \emph{control data} for $f$.
\end{defn}
We now suppose that each of the functions $f\in\F$ has controlled
derivatives. It follows that, in the notation of the previous section, the functions in $\t{\F}$ (and hence in
$R_n$) also have controlled derivatives. Note that, because of the presence of $\exp$, this
assumption holds for the examples of LPB structures given in Section \ref{lpbstructures}.
\begin{propn}
Suppose that $f:U\to M$ is implicitly defined over $\F$. Then $f$ has
controlled derivatives.
\end{propn}
\begin{proof}
Let $g_1,\ldots,g_m\in R_{n+m}$ and $\phi_1,\ldots,\phi_m:U\to M$ witness
the fact that $f$ is implicitly defined. Since $g_1,\ldots,g_m$ have
controlled derivatives, there is a continuous definable function
$\omega:M^n\to M$ and $C_i\in M,E_i\in \N$ such that for each
$i=1,\ldots,m$ and all $\alpha \in N^n$,
$$
|D^{\alpha}g_i(\xbar,\ybar)|\le C_{|\alpha|}\cdot
\omega(\xbar,\ybar)^{E_{|\alpha|}} \textrm{ for all } \la
\xbar,\ybar\ra \in M^{n+m}.
$$
Let $\Delta$ be the determinant of the matrix
$$
\left( \begin{array}{ccc}
\d{g_1}{y_1} &\ldots & \d{g_1}{y_m} \\
\vdots & & \vdots \\
\d{g_m}{y_1} &\ldots & \d{g_m}{y_m} 
\end{array}\right).
$$
We will show by induction on $|\alpha|$ that there are $C'_{|\alpha|}\in
M,E'_{|\alpha|}\in \N$ such that for each $i$ and all $\xbar \in U$,
$$
|D^{\alpha}\phi_i(\xbar)|\le C'_{|\alpha|}\left(
  \frac{\omega(\xbar,\phi_1(\xbar),\ldots,\phi_m(\xbar))}{\Delta(\xbar,\phi_1(\xbar),\ldots,\phi_m(\xbar))}\right)^{E'_{|\alpha|}},
$$ 
which suffices as $f$ is one of the $\phi_i$.

Suppose first that $|\alpha|=1$. We write $\bar{\phi}(\xbar):=\la
\phi_1(\xbar),\ldots,\phi_m(\xbar)\ra$. Since the derivative
$\d{\phi_i}{y_j}(\xbar)$ has the form
$$
\frac{\textrm{polynomial in }\d{g_l}{y_k}\textrm{ evaluated at }\la
\xbar,\bar{\phi}(\xbar)\ra,\textrm{ for various
}k,l}{\Delta(\xbar,\bar{\phi}(\xbar))},
$$
the required $C'_1,E'_1$ clearly exist.

Now suppose that $|\alpha|>1$. By the chain rule,
$D^{\alpha}\phi_i(\xbar)$ has the form
$$
\frac{
\begin{array}{l}
\textrm{polynomial in }D^{\beta}g_j\textrm{ evaluated at }\la
  \xbar,\bar{\phi}(\xbar)\ra\textrm{ and
  }D^{\beta'}\phi_k(\xbar),\\
\textrm{ for various
  }j,k,\beta,\beta'\textrm{ with }|\beta|\le
  |\alpha|,|\beta'|<|\alpha|
\end{array}
}{\Delta(\xbar,\bar{\phi}(\xbar))^d},
$$
and by the induction hypothesis, we can find suitable
$C'_{|\alpha|},E'_{|\alpha|}$.
\end{proof}
Combining this with Corollary \ref{piecewiseimplicit}, we obtain
\begin{cor}\label{controlledderivatives}
Suppose that $f:U\to M$ is a smooth definable function. Then there are
definable open sets $U_1,\ldots,U_k\subseteq U$ with $\dim (U\setminus
\bigcup U_i)<n$ such that for each $i=1,\ldots,k$, $f|_{U_i}$ has
controlled derivatives.
\end{cor}
\begin{rmk} In polynomially bounded structures, all smooth
  functions have controlled derivatives. It seems feasible that a more
  careful analysis of the derivatives of implicit functions may show
  that exponents of the form $|\alpha|$ are preserved. This could lead
  to new results in the polynomially bounded case.
\end{rmk}

\end{document}